\documentclass[a4paper, 11pt]{article}
\usepackage{mathrsfs}
\usepackage{mathrsfs}
\usepackage{amsmath,amssymb}
\usepackage{amsfonts}
\usepackage[T1]{fontenc}
\usepackage[latin9]{inputenc}
\usepackage{amsthm}
\usepackage{graphicx}
\usepackage{subcaption}
\usepackage{amsmath}

\usepackage{amsthm}
\usepackage{array}
\usepackage{cases}
\usepackage{enumerate,enumitem,todonotes}
\usepackage{thmtools}
\usepackage{thm-restate}
\usepackage{setspace}  
\usepackage{changes}
\usepackage{float} 

\makeatletter
\def\th@plain{%
  \itshape 
}
\makeatother

\makeatletter
\renewenvironment{proof}[1][\proofname]{\par
  \pushQED{\qed}%
  \normalfont \topsep6\p@\@plus6\p@\relax
  \trivlist
  \item[\hskip\labelsep
        \bfseries
    #1\@addpunct{.}]\ignorespaces
}{%
  \popQED\endtrivlist\@endpefalse
}
\makeatother

\newtheorem{theorem}{Theorem}[section]

\numberwithin{equation}{section}

\newtheorem{thm}{Theorem}[section]
\newtheorem{cor}[thm]{Corollary}

\newtheorem{lemma}[thm]{Lemma}
\newtheorem{example}[thm]{Example}

\newtheorem{pblm}[thm]{Problem}

\numberwithin{equation}{section}

\usepackage[varg]{txfonts}
\usepackage{graphicx, epsfig}
\usepackage{enumerate}
\usepackage[square, numbers, sort&compress]{natbib}
\ifx\pdfoutput\undefined
 \usepackage[dvipdfm,%
  pdfstartview=FitH,
  bookmarks=true,%
  bookmarksnumbered=true,
  bookmarksopen=true,
  plainpages=false,%
  pdfpagelabels,%
  colorlinks=true,
  linkcolor=blue,
  citecolor=blue,%
  urlcolor=black,
  pdfborder=001]{hyperref}
  \AtBeginDvi{}  
\else
 \usepackage[pdftex,%
  pdfstartview=FitH,
  bookmarks=true,%
  bookmarksnumbered=true,
  bookmarksopen=true,
  plainpages=false,%
  pdfpagelabels,%
  colorlinks=true,
  linkcolor=blue,
  citecolor=blue,%
  urlcolor=black,
  pdfborder=001]{hyperref}%
\fi

\usepackage[left]{lineno}
\RequirePackage[normalem]{ulem} 
\RequirePackage{color}\definecolor{RED}{rgb}{1,0,0}\definecolor{BLUE}{rgb}{0,0,1} 
\providecommand{\DIFaddbegin}{} 
\providecommand{\DIFaddend}{} 
\providecommand{\DIFdelbegin}{} 
\providecommand{\DIFdelend}{} 
\providecommand{\DIFaddbeginFL}{} 
\providecommand{\DIFaddendFL}{} 
\providecommand{\DIFdelbeginFL}{} 
\providecommand{\DIFdelendFL}{} 
\newcommand{\DIFscaledelfig}{0.5}
\RequirePackage{settobox} 
\RequirePackage{letltxmacro} 
\newsavebox{\DIFdelgraphicsbox} 
\newlength{\DIFdelgraphicswidth} 
\newlength{\DIFdelgraphicsheight} 
\LetLtxMacro{\DIFOincludegraphics}{\includegraphics} 
\newcommand{\DIFaddincludegraphics}[2][]{{\color{blue}\fbox{\DIFOincludegraphics[#1]{#2}}}} 
\newcommand{\DIFdelincludegraphics}[2][]{
\sbox{\DIFdelgraphicsbox}{\DIFOincludegraphics[#1]{#2}}
\settoboxwidth{\DIFdelgraphicswidth}{\DIFdelgraphicsbox} 
\settoboxtotalheight{\DIFdelgraphicsheight}{\DIFdelgraphicsbox} 
\scalebox{\DIFscaledelfig}{
\parbox[b]{\DIFdelgraphicswidth}{\usebox{\DIFdelgraphicsbox}\\[-\baselineskip] \rule{\DIFdelgraphicswidth}{0em}}\llap{\resizebox{\DIFdelgraphicswidth}{\DIFdelgraphicsheight}{
\setlength{\unitlength}{\DIFdelgraphicswidth}
\begin{picture}(1,1)
\thicklines\linethickness{2pt} 
{\color[rgb]{1,0,0}\put(0,0){\framebox(1,1){}}}
{\color[rgb]{1,0,0}\put(0,0){\line( 1,1){1}}}
{\color[rgb]{1,0,0}\put(0,1){\line(1,-1){1}}}
\end{picture}
}\hspace*{3pt}}} 
} 
\LetLtxMacro{\DIFOaddbegin}{\DIFaddbegin} 
\LetLtxMacro{\DIFOaddend}{\DIFaddend} 
\LetLtxMacro{\DIFOdelbegin}{\DIFdelbegin} 
\LetLtxMacro{\DIFOdelend}{\DIFdelend} 
\DeclareRobustCommand{\DIFaddbegin}{\DIFOaddbegin \let\includegraphics\DIFaddincludegraphics} 
\DeclareRobustCommand{\DIFaddend}{\DIFOaddend \let\includegraphics\DIFOincludegraphics} 
\DeclareRobustCommand{\DIFdelbegin}{\DIFOdelbegin \let\includegraphics\DIFdelincludegraphics} 
\DeclareRobustCommand{\DIFdelend}{\DIFOaddend \let\includegraphics\DIFOincludegraphics} 
\LetLtxMacro{\DIFOaddbeginFL}{\DIFaddbeginFL} 
\LetLtxMacro{\DIFOaddendFL}{\DIFaddendFL} 
\LetLtxMacro{\DIFOdelbeginFL}{\DIFdelbeginFL} 
\LetLtxMacro{\DIFOdelendFL}{\DIFdelendFL} 
\DeclareRobustCommand{\DIFaddbeginFL}{\DIFOaddbeginFL \let\includegraphics\DIFaddincludegraphics} 
\DeclareRobustCommand{\DIFaddendFL}{\DIFOaddendFL \let\includegraphics\DIFOincludegraphics} 
\DeclareRobustCommand{\DIFdelbeginFL}{\DIFOdelbeginFL \let\includegraphics\DIFdelincludegraphics} 
\DeclareRobustCommand{\DIFdelendFL}{\DIFOaddendFL \let\includegraphics\DIFOincludegraphics} 
\title{\LARGE Between proper and strong edge-colorings of graphs with maximum degree at most four}
\author{
Aistis Atminas\thanks{Department of Applied Mathematics, School of Mathematics and Physics, Xi'an Jiaotong-Liverpool University, Suzhou, Jiangsu Province, 215123, China, \texttt{aistis.atminas@xjtlu.edu.cn}.}
\and
Xujun Liu\thanks{Department of Applied Mathematics, School of Mathematics and Physics, Xi'an Jiaotong-Liverpool University, Suzhou, Jiangsu Province, 215123, China, \texttt{xujun.liu@xjtlu.edu.cn}; the research of X. Liu was supported by the National Natural Science Foundation of China under grants No.~12401466 and No.~12671413, Jiangsu Young Scientific and Technological Talents Promotion Programme, and Suzhou Young Scientific and Technological Talents Promotion Programme.}
\and
Zihui Xu\thanks{School of Mathematics and Physics, Xi'an Jiaotong-Liverpool University, Suzhou, Jiangsu Province, 215123, China, \texttt{zihui.xu2503@student.xjtlu.edu.cn}.}
}
\date{}

\begin{document}

\maketitle

\begin{abstract}
A $(1^j,2^k)$-packing edge-coloring of a graph $G$ is a partition of its edge set $E(G)$ into $j$ matchings and $k$ induced matchings. It can be viewed as intermediate colorings between proper and strong edge-coloring. The study of this topic on subcubic graphs has drawn considerable attention recently, with many conjectures and open problems proposed and resolved. However, very few results are known for graphs with maximum degree four or more. Henderson, Kennedy, and Santana showed that every graph with degree at most four has a $(1,2^{19})$-packing edge-coloring and a $(1^2,2^{17})$-packing edge-coloring.

In this paper, we improve the results of Henderson, Kennedy, and Santana by showing every graph with degree at most four has a $(1,2^{17})$-packing edge-coloring. As a corollary of our result, every graph with degree at most four has a $(1^2,2^{16})$-packing edge-coloring. 
\end{abstract}

\section{Introduction}

Throughout the paper, all graphs are finite, undirected, and simple unless explicitly
stated otherwise.  We write $\Delta(G)$ for the maximum degree of a graph
$G$, and let $L(G)$ denote its line graph, whose vertices are the edges of
$G$, with two vertices adjacent precisely when the corresponding edges
share an endpoint.  The distance between two edges $e,f\in E(G)$ is the
distance between the corresponding vertices in $L(G)$, and is denoted by
$d_G(e,f)$. Given a non-decreasing sequence
$S=(s_1,\ldots,s_t)$ of positive integers, an
\emph{$S$-packing edge-coloring} of $G$ is a partition
\[
    E(G)=E_1\mathbin{\cup}\cdots\mathbin{\cup}E_t
\]
such that for every pair of distinct edges $e,f \in E_i$ the distance $d_G(e,f) \geq s_i+1$ for 
each $i\in\{1,\ldots,t\}$. We use exponents to indicate repeated entries
in $S$; for example, $(1^j,2^k)$ denotes the sequence consisting of $j$
entries equal to $1$ and $k$ entries equal to $2$.  A matching $M$ is
\emph{induced} if the subgraph induced by the endpoints of the edges in
$M$ has edge set exactly $M$.  Consequently, a
$(1^j,2^k)$-packing edge-coloring is precisely a partition of $E(G)$ into
$j$ matchings and $k$ induced matchings.

The two extreme cases recover two classical notions.  A
$(1^j)$-packing edge-coloring is a proper $j$-edge-coloring, while a
$(2^k)$-packing edge-coloring is a strong $k$-edge-coloring.  Thus,
$(1^j,2^k)$-packing edge-colorings form a natural interpolation between
proper and strong edge-coloring.  The least numbers of colors required in
proper and strong edge-colorings of $G$ are denoted by $\chi'(G)$ and
$\chi'_s(G)$, respectively; equivalently, $\chi'_s(G)=\chi(L(G)^2)$,
where $L(G)^2$ denotes the square of the line graph.  A graph $G$ is
\emph{class~I} if $\chi'(G)=\Delta(G)$.  Strong edge-coloring was
introduced by Fouquet and Jolivet~\cite{FouquetJolivet1983}.  Erd\H{o}s and
Ne\v{s}et\v{r}il~\cite{ErdosNesetril1989} conjectured that every graph of
maximum degree $\Delta$ satisfies
\[
  \chi'_s(G)\leq
  \begin{cases}
    \dfrac{5}{4}\Delta^2,
      & \text{if $\Delta$ is even},\\[4pt]
    \dfrac{5}{4}\Delta^2-\dfrac{1}{2}\Delta+\dfrac{1}{4},
      & \text{if $\Delta$ is odd}.
  \end{cases}
\]
For $\Delta=3$, the conjectured bound $10$ was proved independently by
Andersen~\cite{Andersen1992} and by Hor\'{a}k, Qing, and
Trotter~\cite{HorakQingTrotter1993}.

The systematic study of $S$-packing edge-colorings was initiated by
Gastineau and Togni~\cite{GastineauTogni2019}.  Building on the fact that
every subcubic graph admits a
$(1^3,2)$-packing edge-coloring~\cite{Payan1977,FouquetVanherpe2013}, they
investigated how matching colors can be traded for induced-matching colors.
In particular, they proved that every cubic graph with a $2$-factor is
$(1^2,2^5)$-packing edge-colorable, with the improved bound
$(1^2,2^4)$ for $3$-edge-colorable cubic graphs.  They conjectured that
every subcubic graph is $(1^2,2^4)$-packing edge-colorable and asked
whether every cubic graph is $(1,2^7)$-packing edge-colorable.  Hocquard,
Lajou, and Lu\v{z}ar~\cite{HocquardLajouLuzar2022} came within one
induced-matching color of both bounds: every subcubic graph is
$(1^2,2^5)$- and $(1,2^8)$-packing edge-colorable, and the corresponding
bounds for class~I graphs are $(1^2,2^4)$ and $(1,2^7)$.  They further
conjectured that the $(1,2^7)$ bound holds for every subcubic graph.
Subsequently, Liu, Santana, and Short~\cite{LiuSantanaShort2023} proved the
sharp result that every subcubic multigraph is
$(1,2^7)$-packing edge-colorable.  In the two-matching setting, Liu and
Yu~\cite{LiuYu2024} established the conjectured $(1^2,2^4)$ bound for
every connected subcubic graph with more than $70$ vertices.  The numbers
of induced-matching colors in these two targets cannot be decreased in
general: subdividing one edge of $K_{3,3}$ produces a subcubic graph that
is neither $(1,2^6)$- nor $(1^2,2^3)$-packing
edge-colorable~\cite{HocquardLajouLuzar2022,LiuYu2024}.

By contrast, considerably less is known about mixed packing edge-colorings
for graphs of maximum degree four.  In this case, the
Erd\H{o}s--Ne\v{s}et\v{r}il conjecture predicts $\chi'_s(G)\leq 20$.
Successive general upper bounds of $23$, $22$, and $21$ were obtained by
Hor\'{a}k~\cite{Horak1990}, Cranston~\cite{Cranston2006}, and Huang,
Santana, and Yu~\cite{HuangSantanaYu2018}, respectively.  The last result
immediately yields $(1,2^{20})$- and $(1^2,2^{19})$-packing edge-colorings
simply by viewing one or two induced-matching color classes as matching
color classes.  Henderson, Kennedy, and
Santana~\cite{HendersonKennedySantana2025} then proved that every graph
$G$ with $\Delta(G)\leq 4$ admits a $(1,2^{19})$-packing edge-coloring and
a $(1^2,2^{17})$-packing edge-coloring.

In this paper, we improve both bounds of Henderson, Kennedy, and Santana.
Our results are as follows.

\begin{theorem}\label{maintheorem}
    Every graph with maximum degree $\Delta\le4$ is $(1,2^{17})$-packing edge-colorable.
\end{theorem}

Since an induced matching is a matching, we also obtain another improved bound.

\begin{cor}\label{firstcor}
  Every graph with maximum degree $\Delta\le4$ is $(1^2,2^{16})$-packing edge-colorable.  
\end{cor}


Let $k_i$ be the minimum integer $k_i'$ such that every graph with maximum degree at most four has a $(1^i,2^{k_i'})$-packing edge-coloring, where $i\in\{1,2\}$. Theorem~\ref{maintheorem} proves that $k_1\le17$ and Corollary~\ref{firstcor} shows $k_2 \le 16$. The example below shows that $k_1\ge15$ and $k_2\ge10$.

\begin{figure}
\vspace{-5mm}
\centering
\begin{minipage}{0.27\textwidth}
    \centering
    \includegraphics[width=\linewidth]{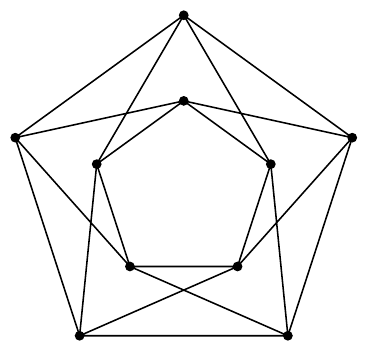}\\
    \textnormal{(a)}
\end{minipage}
\hfill
\begin{minipage}{0.27\textwidth}
    \centering
    \includegraphics[width=\linewidth]{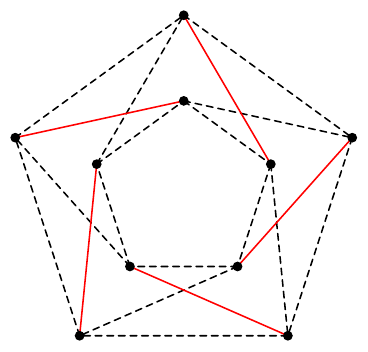}\\
    \textnormal{(b)}
\end{minipage}
\hfill
\begin{minipage}{0.27\textwidth}
    \centering
    \includegraphics[width=\linewidth]{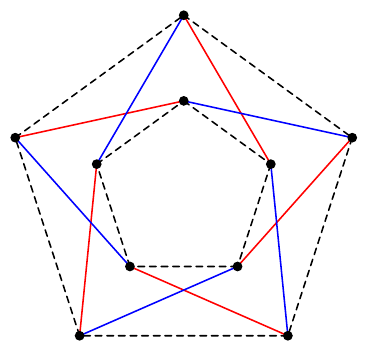}\\
    \textnormal{(c)}
\end{minipage}
\vspace{-3mm}
\caption{A graph with no $(1,2^{14})$-packing edge-coloring and no $(1^2,2^{9})$-packing edge-coloring.}
\label{fig:sharp}
\vspace{-3mm}
\end{figure}

\begin{example}\label{example}
    The graph $G_1$ in the Figure~\ref{fig:sharp}(a) is a graph with maximum degree at most four. Note that $G_1$ has a $(1,2^{15})$-packing edge-coloring and $(1^2,2^{10})$-packing edge-coloring (see Figure~\ref{fig:sharp}(b) and \ref{fig:sharp}(c), where each black dashed edge is assigned a distinct $2$-color, while each set of solid edges shown in the same color is assigned a distinct $1$-color). Since $G_1$ has ten vertices, a $1$-color can
    be used at most five times. Since any two edges of $G_1$ are at
    distance at most two, a $2$-color can be used at most once. Therefore, since $G_1$ has twenty edges, $G_1$ has no $(1,2^{14})$-packing edge-coloring and no $(1^2,2^9)$-packing edge-coloring.
\end{example}

For two vertices $u,v\in V(G)$, the distance $\operatorname{dist}_G(u,v)$ is the length of a shortest path between $u$ and $v$ in $G$. The distance between two edges $e$ and $f$ of $G$ is defined as the distance between their corresponding vertices in the line graph $L(G)$. In this graph, we say that an edge $e=uv$ sees a $1$-color if some edge incident with $u$ or $v$ has that color. Similarly, $e$ sees a $2$-color if some edge at distance at most two from $e$ has that color.
\section{Main proof}
In this section, we show that every graph with maximum degree at most four has a $(1,2^{17})$-packing edge-coloring. We call a $(1,2^{17})$-packing edge-coloring {\em good coloring}. We use $1$ to denote the $1$-color and let $\{2_1, 2_2, \ldots, 2_{17}\}$ be the set of $2$-colors.  


We may assume that $G$ is connected, since if $G$ is disconnected then we apply the same argument to each component of $G$. 
Thus, suppose for a contradiction that there exists a graph of maximum degree at most $4$ that admits no good coloring. 
Let \(G\) be a counterexample with the minimum number of vertices and, subject to this, the maximum number of edges. Let $n=|V(G)|$. Since the result holds for graphs on at most six vertices, we have $n>6$.

\begin{lemma}\label{mindegree}
$\delta(G) \ge 3$.
\end{lemma}
\begin{proof}
    We first show that if $G$ has a 1-vertex, then $G$ has a good coloring. Let $u$ be a 1-vertex of $G$, and let $u_1$ be its unique neighbor. Let $G'=G-u$. Then $G'$ has $n-1$ vertices and $\Delta(G')\le4$. By the minimality of $G$, $G'$ has a good coloring $f$. Since there are at most $12$ edges within distance two from edge $uu_1$ and there are in total 18 colors, we can extend $f$ to $G$. 
    
    We next show that if $\delta(G)=2$, $G$ has a good coloring. Let $u$ be a $2$-vertex with neighbours $u_1,u_2$. Let $G'$ be the graph obtained from $G$ by deleting $u$. Then $G'$ has $n-1$ vertices and $\Delta(G')\le4$. By the minimality of $G$, $G'$ has a good coloring $f$. There are at most $16$ edges within distance two from edge $uu_1$. The same holds for $uu_2$. Since there are in total 18 colors, we can extend $f$ to $G$. 
\end{proof}

\begin{lemma}\label{1-edge-cut}
$G$ has no cut-edge.
\end{lemma}
\begin{proof}
    Suppose not. Let $e=vw$ be a cut-edge. Let the two components of $G-e$ be $G_1$ and $G_2$, with $G_1$ containing vertex $v$, and $G_2$ containing vertex $w$. By Lemma~\ref{mindegree}, we can assume that neither $G_1$ nor $G_2$ is a single vertex. Let $G_1'$ be the graph obtained from $G_1$ by adding the vertex $w'$ and the edge $vw'$, and similarly, let $G_2'$ be the graph obtained from $G_2$ by adding the vertex $v'$ and the edge $v'w$. Since $G_1'$ and $G_2'$ both have fewer vertices than $G$ and have maximum degree at most four, by the minimality of $G$, they both admit good colorings, say $f_1$ and $f_2$, respectively. We may assume that each of $vw'$ and $v'w$ receives a $2$-color. Otherwise, we can recolor $vw'$ or $v'w$ with a $2$-color, since there are at most $12$ edges within distance two from $vw'$ in $G_1'$, the same holds for $v'w$ in $G_2'$ and we have in total seventeen 2-colors. We can then permute the colors in such a way that $vw'$ and $v'w$ receive the same color, say $f_1(vw')=f_2(v'w)=2_1$. Let $S_1$ be the set of edges incident with $v$ in $G_1$ and $S_2$ be the set of edges incident with $w$ in $G_2$. Since 
    $S_1\cup S_2\cup\{vw\}$ contains at most $7$ edges and there are seventeen $2$-colors in total, we can further permute the colors in $G_2'$ so that $f_1(S_1)\cap f_2(S_2)\subseteq\{1\}$, while keeping $2_1$ fixed. Therefore, the colorings of $G_1'$ and $G_2'$ together define a good coloring of $G$, contradicting the choice of $G$.
\end{proof}

\begin{lemma}\label{2-edge-cut}
$G$ has no edge-cut of size two.
\end{lemma}

\begin{proof}

Suppose not. Let $F = \{v_1w_1, v_2w_2\}$ be an edge-cut. We have two cases depending on whether the two edges of $F$ share an endpoint.

\textbf{Case 1:} The two edges of $F$ share an
    endpoint, say $v_1=v_2$. Let the two components of $G-F$ be $G_1$ and $G_2$, with $G_1$ containing vertex $v_1$, and $G_2$ containing vertices $w_1$ and $w_2$. By Lemma~\ref{mindegree}, we can assume that $G_1$ is not a single vertex and $G_2$ has at least three vertices. Let $G_1'$ be the graph obtained from $G_1$ by adding the vertex $w_1',w_2'$ and the edge $v_1w_1',v_1w_2'$, and similarly, let $G_2'$ be the graph obtained from $G_2$ by adding the vertex $v_1'$ and the edge $v_1'w_1,v_1'w_2$. Since $G_1'$ and $G_2'$ both have fewer vertices than $G$ and have maximum degree at most four, by the minimality of $G$, they both admit good colorings, say $f_1$ and $f_2$, respectively. We may assume that each of $v_1w_1', v_1w_2', v_1'w_1$ and $v_1'w_2$ receives a $2$-color, since otherwise we can recolor these edges with $2$-colors. To see this, in $G_1'$, we uncolor $v_1w_1'$ and $v_1w_2'$, each of which sees at most eight $2$-colors. Similarly, we uncolor $v_1'w_1$ and $v_1'w_2$ in $G_2'$, and each of them sees at most fifteen $2$-colors. Since there are seventeen $2$-colors in total, we can color each of $v_1w_1', v_1w_2', v_1'w_1$ and $v_1'w_2$ with a $2$-color.

We can then permute the colors in such a way that $v_1w_1'$ and $v_1'w_1$ have the same color, while $v_1w_2'$ and $v_1'w_2$ also have the same color, say $f_1(v_1w_1')=f_2(v_1'w_1)=2_1$ and $f_1(v_1w_2')=f_2(v_1'w_2)=2_2$. Let $S_1$ be the set of edges incident on $v_1$ in $G_1$ and $S_2$ be the set of edges incident on $w_1$ or $w_2$ in $G_2$. Since $S_1\cup S_2$ contains at most $8$ edges and we still have fifteen $2$-colors excluding $2_1$ and $2_2$, we can further permute the colors in $G_2'$ so that $f_1(S_1)\cap f_2(S_2)\subseteq\{1\}$, while keeping $2_1$ and $2_2$ fixed. Therefore, the colorings of $G_1'$ and $G_2'$ together define a good coloring of $G$.


\textbf{Case 2:} The two edges of $F$ are
    disjoint. By Lemma~\ref{1-edge-cut}, we may assume that $G$ has no cut-edge. Then $G-F$ has exactly two components. By swapping the labels of $v_2$ and $w_2$ if necessary, let $G_1$ and $G_2$ be the components of $G-F$ such that $\{v_1,v_2\}\subseteq V(G_1)$ and ${w_1,w_2}\subseteq V(G_2)$. Now let $G_1'$ be obtained from $G_1$ by adding a new vertex $x_1$ adjacent to both $v_1$ and $v_2$, and let $G_2'$ be obtained from $G_2$ by adding a new vertex $x_2$ adjacent to both $w_1$ and $w_2$. Since $G_1'$ and $G_2'$ both have fewer vertices than $G$ and have maximum degree at most four, by the minimality of $G$, they both admit good colorings, say $f_1$ and $f_2$, respectively. 

We may assume that $v_1x_1$ and $v_2x_1$ are assigned $2$-colors in $G_1'$. To see this, we uncolor $v_1x_1$ and $v_2x_1$, each of which sees at most fifteen $2$-colors. Since there are seventeen $2$-colors in total, we can color each of them with a $2$-color. By the same argument, we can assume that $w_1x_2$ and $w_2x_2$ are assigned $2$-colors in $G_2'$. We can then permute the colors so that $v_1x_1$ and $w_1x_2$ receive the same color, and $v_2x_1$ and $w_2x_2$ receive the same color, say $f_1(v_1x_1)=f_2(w_1x_2)=2_1$ and $f_1(v_2x_1)=f_2(w_2x_2)=2_2$. Let $S_1$ be the set of edges incident with $v_1$ or $v_2$ in $G_1$ and $S_2$ be the set of edges incident with $w_1$ or $w_2$ in $G_2$. Since $S_1\cup S_2$ contains at most $12$ edges and we still have fifteen $2$-colors excluding $2_1$ and $2_2$, we can further permute the colors in $G_2'$ so that $f_1(S_1)\cap f_2(S_2)\subseteq\{1\}$, while keeping $2_1$ and $2_2$ fixed. Therefore, the colorings of $G_1'$ and $G_2'$ together define a good coloring of $G$, where $v_1w_1$ is assigned $2_1$ and $v_2w_2$ is assigned $2_2$.
\end{proof}

\begin{lemma}\label{3v3v}
    Two $3$-vertices cannot be adjacent. 
\end{lemma}
\begin{proof}
    Suppose not. Let $u,v$ be 3-vertices and $uv\in E(G)$.If $u$ and $v$ has no common neighbor, we identify $u$ and $v$ to obtain a graph $G'$ with $\Delta(G')\le4$. Since $G'$ has fewer vertices than $G$, by the minimality of $G$, $G'$ has a good coloring, say $f$. Then we extend $f$ to $G$. Since there are at most sixteen edges within distance two from $uv$ and there are $17$ $2$-colors in total, we can color $uv$ with an available $2$-colors, a contradiction. Then we may assume that $u$ and $v$ have at least one common neighbor. We
    identify $u$ and $v$ into a new vertex $w$ and keep only one copy of each resulting parallel edge to obtain a simple graph $G'$ with $\Delta(G')\le4$. Since $G'$ has fewer vertices than $G$, by the minimality of $G$, $G'$ has a good coloring $f$. We color each edge of $G-uv$ with the $f$-color of its corresponding edge in $G'$ (where $ux$ and $vx$ correspond to $wx$), except that, for every common neighbor $x$ of $u$ and $v$, the edge $vx$ is left uncolored, since otherwise it would get the same color as $ux$. Each such $vx$ has at most fifteen edges within distance two, so it can be colored with an available $2$-color. Since there are at most sixteen edges within distance two from $uv$ and there are $17$ $2$-colors in total, we can color $uv$ with an available $2$-colors.
\end{proof}

\begin{lemma}\label{4-regular}
    $G$ is a $4$-regular graph.
\end{lemma}
\begin{proof}
    By Lemma~\ref{mindegree}, $\delta(G)\ge3$. We only need to prove that $G$ has no $3$-vertex. Suppose not, let $S\ne\emptyset$ be the set of 3-vertices in $G$. Since every vertex of $G$ has degree $3$ or $4$, the handshaking lemma implies that $|S|$ is even. By Lemma~\ref{3v3v}, no two vertices in $S$ are adjacent. Thus, we can pair the vertices of $S$ arbitrarily and add one edge between the two vertices in each pair. Let $G^*$ be the graph obtained in this way. Then $G^*$ is a simple $4$-regular graph with $V(G^*)=V(G)$ and $E(G)\subsetneq E(G^*)$. By the maximality of $|E(G)|$, $G^*$ has a good coloring. Since adding edges can only decrease distances between edges, any $(1,2^{17})$-packing edge-coloring of $G^*$, restricted to $E(G)$, is also a $(1,2^{17})$-packing edge-coloring of $G$, contradicting the choice of $G$. 
\end{proof}

From now on, we may assume that $G$ is $4$-regular.
\begin{lemma}\label{maximum_matching}
    If $|V(G)|$ is even, then $G$ has a perfect matching. If $|V(G)|$ is odd, then $G$ has a near-perfect matching. Moreover, in the latter case, for every vertex $v\in V(G)$, the graph $G-v$ has a perfect matching.
\end{lemma}

\begin{proof}
    For a graph $H$, let $c(H)$ denote the number of components of $H$, and let $\operatorname{odd}(H)$ denote the number of odd components of $H$. We use the Berge-Tutte formula to show that $G$ has a perfect matching or a near-perfect matching. Consider a subset $S \subseteq V(G)$ of $G$ and $G-S$ is a union of $c(G-S)$ components. For a component $C$ of $G-S$, let $\partial_G(C)$ denote the set of edges joining $V(C)$ and $V(G)\setminus V(C)$.

    Assume first that $S\neq \emptyset$. By Lemma~\ref{1-edge-cut} and~\ref{2-edge-cut}, $G$ is 3-edge-connected. Thus, We have $|\partial_G(C)|\ge3$. Moreover, by Lemma~\ref{4-regular}, $G$ is 4-regular. Thus, we have $|\partial_G(C)|=4|V(C)|-2|E(C)|$. Thus, $|\partial_G(C)|$ is even. Hence $|\partial_G(C)|\ge4$ for every component $C$ of $G-S$. Therefore, we have $4c(G-S)\le\sum_C|\partial_G(C)|\le4|S|$, and so $c(G-S)\le|S|$. It follows that 
    \begin{equation}\label{oddcomponent}
     \operatorname{odd}(G-S)\le c(G-S)\le |S|.   
    \end{equation}
    Berge-Tutte formula states that the size of maximum matching in $G$ is 
    \begin{equation}\label{bergetutte}
    \frac{1}{2} \min_{S \subseteq V} (|S|-odd(G-S)+|V|).   
    \end{equation}
    By~\eqref{oddcomponent}, the minimum of~\eqref{bergetutte} is attained at $S=\emptyset$. When $|V(G)|$ is even, the size of maximum matching is $\frac12|V(G)|$, and hence $G$ has a perfect matching. When $|V(G)|$ is odd, the size of maximum matching is $\frac12(|V(G)|-1)$, and hence $G$ has a near-perfect matching.

    Moreover, when $S\ne\emptyset$ and $|V(G)|$ is odd, the equality $|S|=\operatorname{odd}(G-S)$ cannot occur. To see this, if the equality holds, then by~\eqref{oddcomponent}, $|S|=\operatorname{odd}(G-S)=c(G-S)$. Thus all components of $G-S$ are odd components. Let these components be $C_1,C_2,\cdots,C_{|S|}$. Since each $|V(C_i)|$ is odd, $$|V(G)|=|S|+\sum_{i=1}^{|S|}|V(C_i)|=\sum_{i=1}^{|S|}(|V(C_i)|+1)$$ is even. This is a contradiction. Hence, if $|V(G)|$ is odd, then for $S \neq \emptyset$ 
    \begin{equation}\label{nonemptyS}
    |S|-\operatorname{odd}(G-S) \ge 1.    
    \end{equation}
    Now we consider $H=G-v$ for some vertex $v\in V(G)$. For any $S\subseteq V(H)$, by~\eqref{nonemptyS}, $$|S|-\operatorname{odd}(H-S)=|S|-\operatorname{odd}(G-(S\cup\{v\}))=|S\cup\{v\}|-\operatorname{odd}(G-(S\cup\{v\}))-1\ge1-1=0.$$
    The equality holds when $S=\emptyset$. Since $|V(H)|=|V(G)|-1$, by Berge-Tutte formula, $H$ has a matching of size $\frac12|V(H)|=\frac12(|V(G)|-1)$.
\end{proof}

\begin{lemma}\label{16_2-colors}
    Let $v\in V(G)$, and let $M$ be a perfect matching of
    $G$ or of $G-v$. Then $G-v$ has a $(1,2^{16})$-packing edge-coloring
    with $2$-colors $2_1,\ldots,2_{16}$ whose $1$-color class is
    $M\cap E(G-v)$.
\end{lemma}

\begin{proof}

     Let $N(v)=\{a,b,c,d\}$. Thus, $M$ is a perfect or near-perfect matching of $G$. We color all edges of $M$ with color 1. We then consider the remaining edges in order of their distances from $v$, starting with those farther from $v$. For an edge $e=uw\in E(G)$, let $d_v(e)$ denote the distance from $v$ to the closer endpoint of $e$, that is, $d_v(e)=\min\{\operatorname{dist}_G(u,v),\operatorname{dist}_G(w,v)\}$. We order the edges in $E(G)\setminus M$ in nonincreasing order of $d_v(e)$. Edges with the same value of $d_v(e)$ may be considered in any order.

    We claim that the uncolored edges in $G-v$ can be colored greedily, in the above order, using sixteen 2-colors. First, consider an edge $e=uw\in E(G)\setminus M$ with $d_v(e)\ge 2$, as shown in Figure~\ref{fig:1,2_16a}. Since $d_v(e)\ge 2$, neither $u$ nor $w$ belongs to $\{v,a,b,c,d\}$. Without loss of generality, assume that $\operatorname{dist}_G(u,v)=d_v(e)$. Let $u_1$ be the neighbor of $u$ on a shortest path from $u$ to $v$. Then $\operatorname{dist}_G(u_1,v)=d_v(e)-1$. Denote the other neighbors of $u$ by $u_2,u_3$. Thus, $N(u)=\{u_1,u_2,u_3,w\}$ and let $N(w)=\{w_1,w_2,w_3,u\}$. Since $M$ covers every vertex except possibly $v$, each of $\{u_1,u_2,u_3,w_1,w_2,w_3\}$ is incident with an edge of $M$.
    In particular, $u_1$ is incident with exactly three edges outside $M$, each of which is still uncolored when $e$ is colored. Therefore, at most fifteen edges within distance two from $e$ have already been assigned 2-colors, and we can color $e$ with an available 2-color.
    
    
    Now, suppose that $e'$ is incident with a vertex in $N(v)$, say $e'=aa_1\in E(G-v)\setminus M$ where $a\in N(v)$. Let $N(a)=\{v,a_1,a_2,a_3\}$. If $a_1\in N[v]$, say $a_1=b$, as shown in Figure~\ref{fig:1,2_16b}. Let $N(b)=\{v,a,b_1,b_2\}$. Each vertex in $\{a_2, a_3,b_1,b_2\}$ is incident with an edge of $M$. Thus, at most twelve edges within distance two from $e'$ have already been assigned 2-colors, and we can color $e'$ with an available 2-color. Therefore, we may assume $a_1\notin N[v]$. As shown in Figure~\ref{fig:1,2_16a}, let $N(a_1)=\{a,a_{11},a_{12},a_{13}\}$. Each vertex in $\{a_2,a_3,a_{11},a_{12},a_{13}\}$ is incident with an edge of $M$. Since we only consider the edges in $G-v$, at most fifteen edges within distance two from $e'$ have already been assigned 2-colors, and we can color $e'$ with an available 2-color.
\end{proof}

\begin{figure}[htbp]
    \vspace{-6mm}
    \centering
    \begin{subfigure}[c]{0.5\textwidth}
        \centering
        \includegraphics[width=\textwidth]{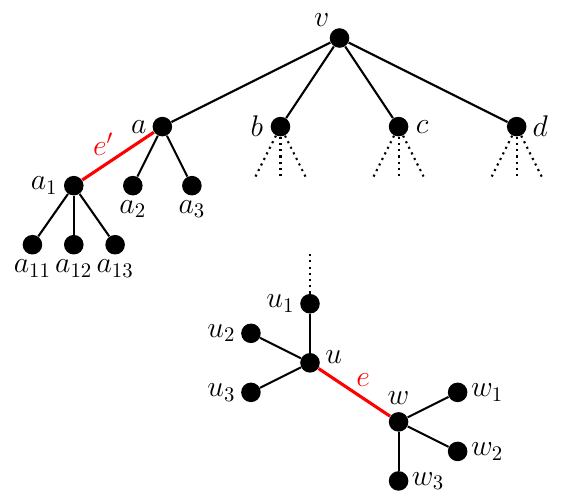}
        \caption{}
        \label{fig:1,2_16a}
    \end{subfigure}
    \hfill
    \begin{subfigure}[c]{0.4\textwidth}
        \centering
        \includegraphics[width=\textwidth]{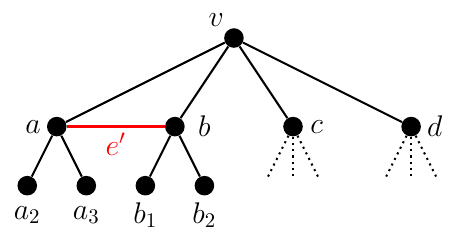}
        \caption{}
        \label{fig:1,2_16b}
    \end{subfigure}
    \vspace{-2mm}
    \caption{Illustration of the edges $e=uw$ and $e'=aa_1$ in the greedy coloring process.}
    \label{fig:1,2_16}
    \vspace{-5mm}
\end{figure}

\begin{lemma}\label{3-cycle}
    $G$ contains no $3$-cycle.
\end{lemma}

\begin{figure}[htbp]
    \vspace{-5mm}
    \centering

    \begin{minipage}[b]{0.42\textwidth}
        \centering
        \includegraphics[width=\textwidth]{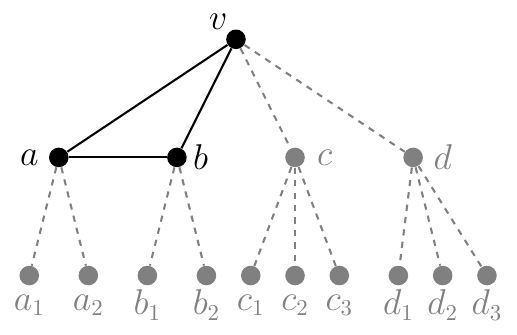}
        \caption{Illustration of Lemma~\ref{3-cycle}}
        \label{fig:3-cycle}
    \end{minipage}
    \hfill
    \begin{minipage}[b]{0.42\textwidth}
        \centering
        \includegraphics[width=\textwidth]{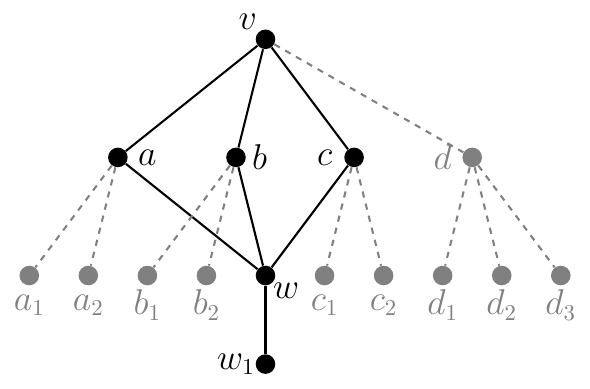}
        \caption{Illustration of Lemma~\ref{K2,3}}
        \label{fig:K2,3}
    \end{minipage}

    \vspace{-5mm}
\end{figure}

\begin{proof}
    Suppose there is a 3-cycle $vabv$ in $G$. Let $N(v)=\{a,b,c,d\}$, $N(a)=\{v,b,a_1,a_2\}$, $N(b)=\{v,a,b_1,b_2\}$, $N(c)=\{v,c_1,c_2,c_3\}$ and $N(d)=\{v,d_1,d_2,d_3\}$, as shown in Figure~\ref{fig:3-cycle}. By Lemma~\ref{maximum_matching}, $G$ or $G-v$ has a perfect matching $M$. We assign the 1-color to every edge of $M$. By Lemma~\ref{16_2-colors}, we use sixteen 2-colors to color the remaining edges in $G-v$.

    If $|V(G)|$ is odd, i.e., $M$ is a perfect matching of $G-v$. Each vertex in $\{a$, $b$, $c$, $d$, $a_1$, $a_2$, $b_1$, $b_2$, $c_1$, $c_2$, $c_3$, $d_1$, $d_2$, $d_3\}$ is incident with an edge of $M$ in $G-v$. For each of $va$ and $vb$, at most twelve edges within distance two of it have already been assigned 2-colors. For each of $vc$ and $vd$, at most fifteen edges within distance two of it have already been assigned 2-colors. We can color the remaining edges in the order $vc,vd,va,vb$.

    If $|V(G)|$ is even, i.e., $M$ is a perfect matching of $G$. First, suppose that one of $va$ and $vb$ is in $M$. By symmetry, we may assume that $va\in M$. Each vertex in $\{b$, $c$, $d$, $b_1$, $b_2$, $c_1$, $c_2$, $c_3$, $d_1$, $d_2$, $d_3\}$ is incident with an edge of $M$ in $G-v$. For $vb$, at most thirteen edges within distance two of it have already been assigned 2-colors. For each of $vc$ and $vd$, at most fifteen edges within distance two of it have already been assigned 2-colors. We can color the remaining edges in the order $vc,vd,vb$.
    Next, suppose that one of $vc$ and $vd$ is in $M$. By symmetry, we may assume that $vc\in M$. Each vertex in $\{a$, $b$, $d$, $a_1$, $a_2$, $b_1$, $b_2$, $d_1$, $d_2$, $d_3\}$ is incident with an edge of $M$ in $G-v$. For each of $va$ and $vb$, at most thirteen edges within distance two of it have already been assigned 2-colors. For each of $vd$, at most sixteen edges within distance two of it have already been assigned 2-colors. We can color the remaining edges in the order $vd$, $va$, and $vb$.
\end{proof}

\begin{lemma}\label{K2,3}
    $G$ contains no $K_{2,3}$.
\end{lemma}
\begin{proof}
    Suppose that $G$ contains a copy of $K_{2,3}$ with bipartition $\{v,w\}$ and $\{a,b,c\}$. Let $N(v)=\{a,b,c,d\}$, $N(a)=\{v,w,a_1,a_2\}$, $N(b)=\{v,w,b_1,b_2\}$, $N(c)=\{v,w,c_1,c_2\}$, $N(d)=\{v,d_1,d_2,d_3\}$, and $N(w)=\{a,b,c,w_1\}$, as shown in Figure~\ref{fig:K2,3}. By Lemma~\ref{3-cycle}, we assume that any two of $a,b,c$ and $d$ are not adjacent. By Lemma~\ref{maximum_matching}, $G$ or $G-v$ has a perfect matching $M$. We assign the 1-color to every edge of $M$. By Lemma~\ref{16_2-colors}, we use sixteen 2-colors to color the remaining edges in $G-v$.

     If $|V(G)|$ is odd, i.e., $M$ is a perfect matching of $G-v$. Each vertex in $\{a$, $b$, $c$, $d$, $w$, $a_1$, $a_2$, $b_1$, $b_2$, $c_1$, $c_2$, $d_1$, $d_2$, $d_3\}$ is incident with an edge of $M$ in $G-v$. Note that at least two of $wa,wb$ and $wc$ are not in $M$. Without loss of generality, assume that $wb,wc\notin M$. The vertices $a_1,a_2,b,c,d$ are incident with at most $3,3,2,2,2$ edges that have already been assigned 2-colors. Since $wb,wc\notin M$, at most one of $wa$ or $ww_1$ has already been assigned a 2-color. Therefore, at most $13$ edges within distance two of $va$ have already been assigned 2-colors.
     Similarly, at most $14$, $14$, and $15$ edges within distance two of $vb$, $vc$, and $vd$, respectively, have already been assigned 2-colors. We can color the remaining edges in the order $vd,vc,vb,va$.

     If $|V(G)|$ is even, i.e., $M$ is a perfect matching of $G$. Note that each vertex in $\{a$, $b$, $c$, $d$, $w$, $a_1$, $a_2$, $b_1$, $b_2$, $c_1$, $c_2$, $d_1$, $d_2$, $d_3\}$ is incident with an edge of $M$ in $G-v$. First, suppose that one of $va,vb$ and $vc$ is in $M$. By symmetry, we may assume that $va\in M$. Note that at least one of $wb$ and $wc$ is not in $M$. Without loss of generality, assume that $wc\notin M$. The vertices $a,c,d,b_1,b_2$ are incident with at most $3,2,2,3,3$ edges that have already been assigned 2-colors. Since $wa, wc\notin M$, at most one of $wb$ or $ww_1$ has already been assigned a 2-color. Therefore, at most $14$ edges within distance two of $vb$ have already been assigned $2$-colors.
     Similarly, at most $15$, and $16$ edges within distance two of $vc$, and $vd$, respectively, have already been assigned 2-colors. We can color the remaining edges in the order $vd,vc,vb$.
     Next, suppose $vd$ is in $M$. Note that at least two of $wa,wb$ and $wc$ are not in $M$. Without loss of generality, assume that $wb,wc\notin M$. The vertices $b,c,d,a_1,a_2$ are incident with at most $2,2,3,3,3$ edges that have already been assigned 2-colors. Since $wb, wc\notin M$, at most one of $wa$ or $ww_1$ has already been assigned a 2-color. Therefore, at most $14$ edges within distance two of $va$ have already been assigned $2$-colors.
     Similarly, at most $15$, and $15$ edges within distance two of $vb$ and $vc$, respectively, have already been assigned 2-colors. We color the remaining edges in the order $vc,vb,va$.
\end{proof}

\begin{figure}[htbp]
    \vspace{-6mm}
    \centering

    \begin{minipage}[t]{0.5\textwidth}
        \centering
        \includegraphics[width=0.9\textwidth]{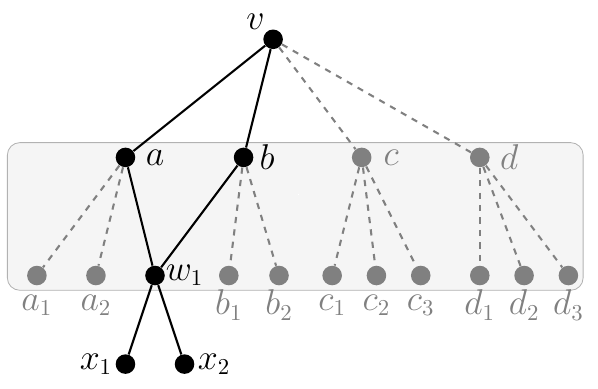}
        \caption{Illustration of Lemma~\ref{C4K2K2}. The gray box contains a pair of edges that are not in $M$ and are at distance at least $3$ from each other.}
        \label{fig:C4K2K2}
    \end{minipage}
    \hfill
    \vspace{-3mm}
    \begin{minipage}[t]{0.40\textwidth}
        \centering
        \includegraphics[width=\textwidth]{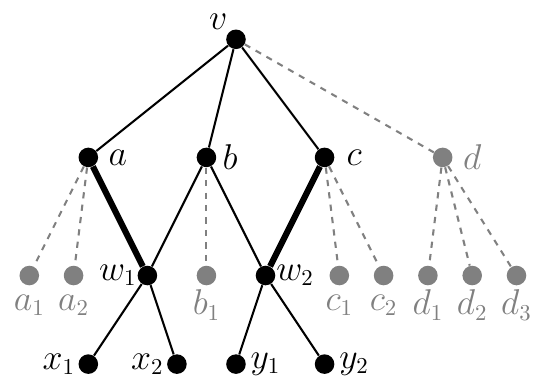}
        \caption{Illustration of Lemma~\ref{C4|C4}. $w_1a,w_2c\notin M$.}
        \label{fig:C4C4}
    \end{minipage}

    \vspace{-2mm}
\end{figure}

Next, we prove $G$ contains no $4$-cycle, after showing a few lemmas. Let $v$ be a vertex of $G$ with $N(v)=\{a,b,c,d\}$. Suppose that there exists a vertex $w_1\ne v$ such that $w_1$ is adjacent to both $a$ and $b$ to form a $4$-cycle $vaw_1bv$. Let $M$ be a perfect matching in $G$ or a perfect matching in $G-v$ guaranteed by Lemma~\ref{maximum_matching}. Let $S$ be the set of edges in $G-v$ that are incident with a vertex of $\{a,b,c,d\}$, see boxed edges in Figure~\ref{fig:C4K2K2}. 

\begin{lemma}\label{C4K2K2}
If there are two edges in $S-M$ that are at distance at least $3$ from each other, then $G$ has a good coloring.
\end{lemma}

\begin{proof}
    Let $N(a)=\{v,w_1,a_1,a_2\}$, $N(b)=\{v,w_1,b_1,b_2\}$, $N(c)=\{v,c_1,c_2,c_3\}$, $N(d)=\{v,d_1,d_2,d_3\}$ and $N(w_1)=\{a,b,x_1,x_2$\}, as shown in Figure~\ref{fig:C4K2K2}. By Lemma~\ref{3-cycle}, we assume that any two of $a,b,c$ and $d$ are not adjacent. We assign 1-color to every edge of $M$. By Lemma~\ref{16_2-colors}, we use sixteen 2-colors, say $\{2_1, 2_2, \ldots, 2_{16}\}$, to color the remaining edges in $G-v$. Since $S-M$ contains two edges at distance $3$ in $G$, we can recolor these two edges with color $2_{17}$.

    If $|V(G)|$ is odd, i.e., $M$ is a perfect matching of $G-v$. Each vertex in $\{a$, $b$, $c$, $d$, $w_1$, $a_1$, $a_2$, $b_1$, $b_2$, $c_1$, $c_2$, $c_3$, $d_1$, $d_2$, $d_3\}$ is incident with an edge of $M$ in $G-v$. Note that at least one of $w_1a$ and $w_1b$ is not in $M$. Without loss of generality, assume that $w_1b\notin M$. The vertices $b,c,d,a_1,a_2$ are incident with at most $2,2,2,3,3$ edges that have already been assigned $2$-colors. Since $w_1b\notin M$, at most two of $w_1a, w_1x_1$, and $w_1x_2$ are assigned $2$-colors. Moreover, $2_{17}$ appears twice in $S$. Therefore, at most $13$ distinct $2$-colors appear within distance two of $va$. Similarly, at most $14$, $14$, and $14$ distinct $2$-colors appear within distance two of $vb$, $vc$, and $vd$, respectively. We can color the remaining edges in the order $vb$, $vc$, $vd$, and $va$.

    If $|V(G)|$ is even, i.e., $M$ is a perfect matching of $G$. Note that each vertex in $\{a$, $b$, $c$, $d$, $w_1$, $a_1$, $a_2$, $b_1$, $b_2$, $c_1$, $c_2$, $c_3$, $d_1$, $d_2$, $d_3\}$ is incident with an edge of $M$ in $G$. First, suppose that one of $va$ and $vb$ is in $M$, say $va\in M$. The vertices $a,c,d,b_1,b_2$ are incident with at most $3,2,2,3,3$ edges that have been assigned $2$-colors. Since $w_1a\notin M$, at most two of $w_1b,w_1x_1$ and $w_1x_2$ are assigned $2$-colors. Therefore, at most $14$ distinct $2$-colors appear within distance two from $vb$. Similarly, at most $15$ and $15$ distinct $2$-colors appear within distance two of $vc$ and $vd$, respectively. We can color the remaining edges in the order $vc$, $vd$, and $vb$.
    Next, suppose that one of $vc$ and $vd$ is in $M$, say $vc\in M$. Note that at least one of $w_1a$ and $w_1b$ is not in $M$, say $w_1b\notin M$. The vertices $b,c,d,a_1,a_2$ are incident with at most $2,3,2,3,3$ edges that have been assigned $2$-colors. Since $w_1b\notin M$, at most two of $w_1a,w_1x_1$ and $w_1x_2$ are assigned $2$-colors. Therefore, at most $14$ distinct $2$-colors appear within distance two from $va$. Similarly, at most $15$ and $15$ distinct $2$-colors appear within distance two of $vb$ and $vd$, respectively. We color the remaining edges in the order $vb$, $vd$, and $va$.
\end{proof}

\begin{lemma}\label{C4|C4}
    Note that at least one of $w_1a$ and $w_1b$ are not in $M$, say $w_1a\notin M$. If there exists a vertex $w_2$ such that $w_2b, w_2c\in E(G)$ and $w_2c\notin M$, then $G$ has a good coloring.
\end{lemma}

\begin{proof}
    Let $N(a)=\{v,w_1,a_1,a_2\}$, $N(b)=\{v,w_1,w_2,b_1\}$, $N(c)=\{v,w_2,c_1,c_2\}$, $N(d)=\{v,d_1,d_2,d_3\}$, $N(w_1)=\{a,b,x_1,x_2\}$ and $N(w_2)=\{b,c,y_1,y_2\}$, as shown in Figure~\ref{fig:C4C4}. By Lemma~\ref{3-cycle}, we assume that any two of $a,b,c$ and $d$ are not adjacent. By Lemma~\ref{maximum_matching}, $G$ or $G-v$ has a perfect matching $M$. We assign the 1-color to every edge of $M$. By Lemma~\ref{16_2-colors}, we use sixteen 2-colors to color the remaining edges in $G-v$.

     If $|V(G)|$ is odd, i.e., $M$ is a perfect matching of $G-v$. Each vertex in $\{a$, $b$, $c$, $d$, $w_1$, $w_2$, $a_1$, $a_2$, $b_1$, $c_1$, $c_2$, $d_1$, $d_2$, $d_3\}$ is incident with an edge of $M$ in $G-v$. Note that at least one of $w_1b$ and $w_2b$ is not in $M$, say $w_1b\notin M$. The vertices $b, c, d, a_1, a_2$ are incident with at most $2,2,2,3,3$ edges that have been assigned $2$-colors. Since $w_1b\notin M$, at most two of $w_1a, w_1x_1$ and $w_1x_2$ are assigned $2$-colors. Therefore, at most $14$ distinct $2$-colors appear within distance two of $va$.
     Similarly, at most $13$, $15$, and $15$ edges within distance two of $vb$, $vc$, and $vd$, respectively, have already been assigned 2-colors. We can color the remaining edges in the order $vc,vd,va,vb$.

     If $|V(G)|$ is even, i.e., $M$ is a perfect matching of $G$. By Lemma~\ref{3-cycle}, $w_1w_2\notin E$. By lemma~\ref{K2,3}, $cw_1, aw_2\notin E(G)$. Thus, the distance between $w_1a$ and $w_2c$ is three, and we can recolor these two edges with color $2_{17}$. Note that each vertex in $\{a$, $b$, $c$, $d$, $w_1$, $w_2$, $a_1$, $a_2$, $b_1$, $c_1$, $c_2$, $d_1$, $d_2$, $d_3\}$ is incident with an edge of $M$ in $G-v$. First, suppose that one of $va$ and $vc$ is in $M$, say $va\in M$. The vertices $a,c,d,b_1$ are incident with at most $3,2,2,3$ edges that have been assigned $2$-colors. Since $w_1a\notin M$, at most two of $w_1b, w_1x_1$ and $w_1x_2$ are assigned $2$-colors. Since $w_2c\notin M$, at most two of $w_2b,w_2y_1$, and $w_2y_2$ are assigned $2$-colors. Moreover, $w_1a$ and $w_2c$ use the same color. Therefore, at most $13$ distinct $2$-colors appear within distance two of $vb$. Similarly, at most $15$ and $15$ distinct $2$-colors appear within distance two of $vc$ and $vd$, respectively. We can color the remaining edges in the order $vc,vd,vb$. 
     Next, suppose that $vb$ is in $M$. Since $w_1a$ and $w_2c$ use the same color, at most $14$, $14$, and $15$ distinct $2$-colors appear within distance two of $va$, $vc$, and $vd$, respectively. We can color the remaining edges in the order $vd,va,vc$. 
     Finally, suppose that $vd$ is in $M$. Note that at least one of $w_1b$ and $w_2b$ is not in $M$ say $w_1b\notin M$. The vertices $b,c,d,a_1,a_2$ are incident with at most $2,2,3,3,3$ edges that have been assigned $2$-colors. Since $w_1b\notin M$, at most two of $w_1a, w_1x_1$ and $w_1x_2$ are assigned $2$-colors. Moreover, $w_1a$ and $w_2c$ use the same color. Therefore, at most $14$ distinct $2$-colors appear within distance two of $va$. Similarly, at most $13$, and $15$ distinct $2$-colors appear within distance two of $vb$, and $vc$, respectively. We can color the remaining edges in the order $vc$, $va$, and $vb$.
\end{proof}

\begin{lemma}\label{4-cycle}
    $G$ contains no 4-cycle.
\end{lemma}

\begin{figure}[htbp]
    \vspace{-3mm}
    \centering

    \begin{minipage}[b]{0.45\textwidth}
        \centering
        \includegraphics[width=\textwidth]{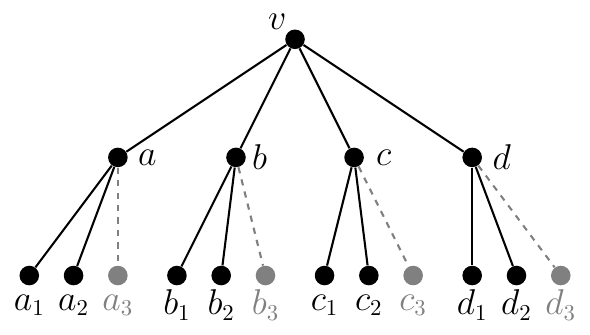}
        \caption{Illustration of Lemma~\ref{4-cycle}}
        \label{fig:C4}
    \end{minipage}
    \hfill
    \begin{minipage}[b]{0.45\textwidth}
        \centering
        \includegraphics[width=\textwidth]{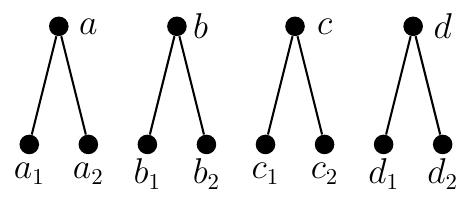}
        \caption{Illustration of Lemma~\ref{4-cycle}.}
        \label{fig:C4_H}
    \end{minipage}

    \vspace{-2mm}
\end{figure}

\begin{proof}

    Suppose $v$ is a vertex of a 4-cycle. Let $N(v)=\{a,b,c,d\}$, $N(a)=\{v,a_1,a_2,a_3\}$, $N(b)=\{v,b_1,b_2,b_3\}$, $N(c)=\{v,c_1,c_2,c_3\}$ and $N(d)=\{v,d_1,d_2,d_3\}$, as shown in Figure~\ref{fig:C4}. By Lemma~\ref{3-cycle}, we assume that any two of $a,b,c$ and $d$ are not adjacent. By Lemma~\ref{maximum_matching}, $G$ or $G-v$ has a perfect matching $M$. We assign the 1-color to every edge of $M$. By Lemma~\ref{16_2-colors}, we use sixteen 2-colors to color the remaining edges in $G-v$.

    Note that each of $a,b,c,d$ is incident with an edge of $M$. Without loss of generality, we may assume that $\{aa_1$, $aa_2$, $bb_1$, $bb_2$, $cc_1$, $cc_2$, $dd_1$, $dd_2\}\notin M$. Let $E_0=\{aa_1$, $aa_2$, $bb_1$, $bb_2$, $cc_1$, $cc_2$, $dd_1$, $dd_2\}$ and $V_0=\{a_1$, $a_2$, $b_1$, $b_2$, $c_1$, $c_2$, $d_1$, $d_2$ $\}$, as shown in Figure~\ref{fig:C4_H}. Let $H$ be the subgraph of $G$ induced by $E_0$. Throughout this proof, adjacency and distances are taken in $G$. By Lemma~\ref{K2,3}, we may assume that $G$ contains no $K_{2,3}$. Thus, each vertex in $V_0$ is adjacent to at most two of $\{a,b,c,d\}$ and any two vertices in $\{a,b,c,d\}$ have at most one common neighbor in $V_0$.
    By Lemma~\ref{C4|C4}, we may assume that $H$ contains no $P_6$. By Lemma~\ref{C4K2K2}, if $H$ contains two edges at distance at least $3$, then we are done. Thus, it remains to show that such two edges always exist.

    Note that the subgraph $H$ contains the four paths $a_1aa_2$, $b_1bb_2$, $c_1cc_2$, and $d_1dd_2$. Here $a,b,c,d$ are pairwise distinct and pairwise nonadjacent. By the previous observation and since $H$ contains no $P_6$, if two of these paths have a common vertex, then they have no common vertex with the other two paths. Therefore, up to symmetry, we need only consider the following three cases.

    \textbf{Case 1:} $H$ consists of four disjoint copies of $P_3$, i.e., vertices in $V_0$ are pairwise distinct. Each of $a,b,c,d$ is incident with exactly one edge of $G-v$ not in $E_0$. Since $|V_0|=8$, at least four vertices in $V_0$ are not adjacent with any of $a,b,c$ and $d$ via an edge not in $E_0$. Without loss of generality, let $a_1$ be such a vertex. Note that $a_1$ is adjacent to at most three vertices in $V_0\setminus\{a_1,a_2\}$, and $a$ is adjacent to at most one vertex in $V_0\setminus\{a_1,a_2\}$, while neither $a$ nor $a_1$ is adjacent to any of $b,c,d$. Since $|V_0\setminus\{a_1,a_2\}|=6$, there exists an edge in $E_0$ at distance at least three from $aa_1$. By Lemma~\ref{C4K2K2}, $G$ has a good coloring.

    \textbf{Case 2:} $H$ consists of one copy of $P_5$ and two copies of $P_3$, say $a_1=b_1$. By Lemma~\ref{C4|C4}, neither $a$ nor $b$ is adjacent to any of $c_1,c_2,d_1,d_2$. Moreover, $c$ is adjacent to at most one in $\{d_1,d_2\}$, since otherwise $G$ would contain a copy of $K_{2,3}$. By symmetry, we may assume that $cd_1\notin E(G)$. Note that $d_1$ is adjacent to at most three vertices in $V_0\setminus\{d_1,d_2\}$, and $d$ is adjacent to at most one vertex in $V_0\setminus\{d_1,d_2\}$, while neither $d$ nor $d_1$ is adjacent to any of $a,b,c$. Since $|V_0\setminus\{d_1,d_2\}|=5$, there exists an edge in $E_0$ at distance at least three from $dd_1$. By Lemma~\ref{C4K2K2}, $G$ has a good coloring.

    \textbf{Case 3:} $H$ consists of two disjoint copies of $P_5$, say $a_1=b_1$ and $c_1=d_1$.By Lemma~\ref{K2,3}, $a$ is not adjacent to $b_2$ and $b$ is not adjacent to $a_2$. By Lemma~\ref{C4|C4}, neither $a$ nor $b$ is adjacent to any of $c_1,c_2,d_2$, and neither $c$ nor $d$ is adjacent to any of $a_1,a_2,b_2$. Note that $a_2$ is adjacent to at most three vertices in $V_0\setminus\{a_1,a_2\}$, while neither $a$ nor $a_2$ is adjacent to any of $b,c,d$. Since $|V_0\setminus\{a_1,a_2\}|=4$, there exists an edge in $E_0$ at distance at least three from $aa_2$. By Lemma~\ref{C4K2K2}, $G$ has a good coloring.
\end{proof}

\begin{lemma}\label{3pairs}
    Let $v$ be a vertex of $G$ with $N(v)=\{a,b,c,d\}$. Let $M$ be a perfect matching in $G$ or $G-v$ guaranteed by Lemma~\ref{maximum_matching}. Let $S$ be the set of edges in $G-v$ that are incident with a vertex of $\{a,b,c,d\}$. If there exist three edges $e_1,e_2,$ and $e_3$ in $S-M$ such that the distance between any two of them is three, then $G$ admits a good coloring.
\end{lemma}

\begin{proof}
    Let $N(v)=\{a,b,c,d\}$. Note that $M$ is a perfect or near-perfect matching of $G$. We assign the 1-color to every edge of $M$. By Lemma~\ref{16_2-colors}, we use sixteen 2-colors to color the remaining edges in $G-v$. We recolor $e_1,e_2$ and $e_3$ with $2_{17}$.

    If $|V(G)|$ is odd, i.e., $M$ is a perfect matching of $G-v$. Each vertex in $\{a$, $b$, $c$, $d$, $a_1$, $a_2$, $a_3$, $b_1$, $b_2$, $b_3$, $c_1$, $c_2$, $c_3$, $d_1$, $d_2$, $d_3\}$ is incident with an edge of $M$ in $G-v$. Therefore, for each edge in $\{va,vb,vc,vd\}$, at most $13$ distinct $2$-colors appear within distance two from it. We can color the remaining edges in the order $va$, $vb$, $vc$, and $vd$.

    If $|V(G)|$ is even, i.e., $M$ is a perfect matching of $G$. Without loss of generality, we may assume $va\in M$. Each vertex in $\{a$, $b$, $c$, $d$, $b_1$, $b_2$, $b_3$, $c_1$, $c_2$, $c_3$, $d_1$, $d_2$, $d_3\}$ is incident with an edge of $M$ in $G-v$. Therefore, for each edge in $\{vb,vc,vd\}$, at most $14$ distinct $2$-colors appear within distance two from it. We can color the remaining edges in the order $vb$, $vc$, and $vd$.
\end{proof}

\begin{proof}[Proof of Theorem~\ref{maintheorem}:]

Let $v$ be a vertex of $G$. Let $N(v)=\{a,b,c,d\}$, $N(a)=\{v,a_1,a_2,a_3\}$, $N(b)=\{v,b_1,b_2,b_3\}$, $N(c)=\{v,c_1,c_2,c_3\}$ and $N(d)=\{v,d_1,d_2,d_3\}$. By Lemma~\ref{3-cycle} and ~\ref{4-cycle}, we may assume that $G$ does not contain 3-cycle and 4-cycle. By Lemma~\ref{maximum_matching}, $G$ or $G-v$ has a perfect matching $M$. We assign the 1-color to every edge of $M$. By Lemma~\ref{16_2-colors}, we use sixteen 2-colors to color the remaining edges in $G-v$.
Note that each vertex in $\{a$, $b$, $c$, $d$, $a_1$, $a_2$, $a_3$, $b_1$, $b_2$, $b_3$, $c_1$, $c_2$, $c_3$, $d_1$, $d_2$, $d_3\}$ is incident with an edge of $M$ in $G-v$. Without loss of generality, we assume that $aa_1,aa_2,bb_1,bb_2,cc_1,cc_2,dd_1,dd_2\notin M$.Let $V_0=\{a_1$, $a_2$, $b_1$, $b_2$, $c_1$, $c_2$, $d_1$, $d_2\}$ and let $H$ be the graph induced by $\{a$, $b$, $c$, $d\}\cup V_0$. Note that each of $a,b,c,$ and $d$ has degree $2$ in $H$, since otherwise $G$ would contain 3-cycle or 4-cycle. By Lemma~\ref{3pairs}, if $H$ contains three edges that pairwise at distance at least three, then $G$ has a good coloring. Thus, it remains to show that $H$ always contains three such edges.

First, suppose that there is a $3^+$-vertex in $H$. Without loss of generality, let $a_1$ be such a vertex and $b_1\in N(a_1)$. Then $a_1b_2\notin E(G)$, since otherwise $G$ would contain a 4-cycle. Without loss of generality, let $c_1$ be adjacent to $a_1$. Note that $b_1c_1\notin E(G)$, since otherwise $G$ would contain a 3-cycle. Moreover, $a_2b_1,a_2c_1\notin E(G)$, since otherwise, $G$ would contain a 4-cycle. Therefore, the distance between any two of the edges $aa_2$, $bb_1$, and $cc_1$ is three. By Lemma~\ref{3pairs}, $G$ has a good coloring.

Next, suppose that there is no $3^+$-vertex in $H$. Note that at least one vertex in $V_0$ is a $2$-vertex in $H$, since otherwise the distance between any two of the edges $aa_1$, $bb_1$, and $cc_1$ is three, and by Lemma~\ref{3pairs}, $G$ has a good coloring. Without loss of generality, we may assume that $a_1$ is a $2$-vertex in $H$. Since $H$ contains no $3$-cycle, $a_1$ is adjacent to exactly one vertex in $\{b_1$, $b_2$, $c_1$, $c_2$, $d_1$, $d_2\}$. Without loss of generality, we assume that $a_1$ is not adjacent to any of $b_2$, $c_1$, $c_2$, $d_1$, $d_2$. Now $b_2$ is adjacent to at most one vertex in $\{a_2, c_1,c_2,d_1,d_2\}$. Without loss of generality, we assume that $b_2$ is not adjacent to $c_1$. Therefore, the distance between any two of the edges $aa_1$, $bb_2$, and $cc_1$ is three. By Lemma~\ref{3pairs}, $G$ has a good coloring.
\end{proof}

\section{Future Directions}

Recall that $k_i$ is the minimum integer $k_i'$ such that every graph with maximum degree at most four has a $(1^i,2^{k_i'})$-packing edge-coloring, where $i\in\{1,2\}$. In this paper, we proved $15 \le k_1 \le 17$ and $10 \le k_2 \le 16$. 

We end the paper by posting an open problem:

\begin{pblm}
What is the exact value of $k_1$ and $k_2$?     
\end{pblm}


\end{document}